\documentclass[twoside,notitlepage,11pt]{article}

\usepackage{amssymb}
\usepackage[leqno]{amsmath}
\usepackage{amsfonts}
\usepackage{amsopn}
\usepackage{amstext}
\usepackage{amsthm}
\usepackage[all]{xy}

\newenvironment{pf}{\begin{proof}}{\end{proof}}

\newcommand{\lam}{{\lambda}}
\newcommand{\al}{\alpha}

\renewcommand{\phi}{\varphi}
\renewcommand{\rho}{\varrho}

\newcommand{\rest}{\restriction}

\newcommand{\loe}{\leqslant}
\newcommand{\goe}{\geqslant}

\newcommand{\subs}{\subseteq}

\newcommand{\cl}{\operatorname{cl}}

\newcommand{\w}{\operatorname{w}}

\newcommand{\conv}{\operatorname{conv}}

\newcommand{\id}{\operatorname{id}}

\newcommand{\liminv}{\varprojlim}

\newcommand{\setof}[2]{\{#1\colon #2\}}

\newcommand{\sett}[2]{\{#1\}_{#2}}
\newcommand{\map}[3]{#1\colon #2 \to #3} 
\newcommand{\img}[2]{#1[#2]} 

\newcommand{\im}{\operatorname{im}}

\newcommand{\norm}[1]{\|#1\|}

\newcommand{\R}{\ensuremath{\mathcal R}}

\newcommand{\invsys}[5]{\langle {#1}_{#4},{#2}_{#4}^{#5},#3 \rangle}

\renewcommand{\S}{\mathbb S}

\newcommand{\alom}{{\al<\omega_1}}

\newcommand{\pli}{Plichko}

\newtheorem{tw}{Theorem}[section]

\newtheorem{lm}[tw]{Lemma}

\newtheorem{twH}{Haydon's}

\theoremstyle{definition}

\theoremstyle{remark}
\newtheorem{uwgi}[tw]{Remark}

\title{Spaces of continuous functions over Dugundji compacta}
\author{
{\sc Taras Banakh}\\ \\
{\small Institute of Mathematics,}\\
{\small Akademia \'Swi\c etokrzyska}\\ {\small (Kielce, Poland)}\\
{\small\it and}\\
{\small Institute of Mathematics,}\\
{\small Ivan Franko Lviv National University}\\ {\small (Lviv, Ukraine)}
\and
{\sc Wies{\l}aw Kubi\'s}
\\ \\
{\small Institute of Mathematics,}\\
{\small Akademia \'Swi\c etokrzyska}\\ {\small (Kielce, Poland)}
}

\begin{document}
\maketitle

\begin{abstract} We show that for every Dugundji compact $K$ of weight $\aleph_1$ the Banach space $C(K)$ is $1$-\pli~and the space $P(K)$ of probability measures on $K$ is Valdivia compact. Combining this result with the existence of a non-Valdivia compact group, we answer a question of Kalenda.

\vspace{3mm}

\end{abstract}

\section{Introduction}

Given an infinite-dimensional Banach space $E$, it is an important question how many nontrivial projections it does have. A {\em projection} is, by definition, a bounded linear operator $\map PEE$ such that $PP=P$. A projection $P$ is {\em nontrivial} if both $\im P$ and $\ker P$ are infinite-dimensional.
There exist {\em indecomposable} Banach spaces, i.e. spaces on which every projection is trivial; by \cite{Koszmider} such spaces can be even of the form $C(K)$ for a suitable (non-metrizable) compact $K$. 

We are interested in non-separable Banach spaces which have ``many" nontrivial norm one projections onto separable subspaces.
Note that every retraction of a compact space $K$ induces a norm one projection on $C(K)$.
Thus, it is natural to ask whether there is a compact space $K$ having few retractions onto metrizable subsets, while its Banach space $C(K)$ has ``many" norm one projections onto separable subspaces.
Of course such a question is not precise. A reasonable yet large enough class of Banach spaces with ``many" norm one projections seems to be the class of {\em $1$-\pli~spaces} (see the definitions below).
On the other hand, one should say that a compact $K$ has {\em few retractions} if for every retraction $\map rKK$ with $X=\img rK$ being second countable, there are some restrictive conditions on the topological type of $X$. For example: the first cohomology group $H^1(X)$ be trivial. This is motivated by \cite{KU}, where a compact (connected Abelian) group $G$ with this property is described. Trying to investigate this particular $C(G)$ space, one can look for a topological property of $G$ imposing the existence of many projections on $C(G)$. It turns out that the property of being {\em Dugundji compact} is good enough. Namely, we show  that $C(K)$ is a $1$-\pli~space for every Dugundji compact $K$ of weight $\loe\aleph_1$. In particular, it follows that the space of probability measures $P(K)$ may be Valdivia compact (the property dual to being $1$-\pli), while at the same time $K$ can be relatively far from Valdivia compacta (again witnessed by the compact group $G$ from \cite{KU}). This answers a question of Ond\v rej Kalenda \cite[Question 5.1.10]{Kalenda} in the negative.

We do not know whether $C(K)$ is $1$-\pli, for a Dugundji compact $K$ of weight $>\aleph_1$. One can prove that in this case $P(K)$ is a retract of a Valdivia compact and $C(K)$ is isomorphic to a $1$-complemented subspace of a \pli\ space. However, it is an open question whether Valdivia compacta are stable under retracts and whether \pli\ spaces are stable under complemented subspaces.
One should mention that every \pli~space admits an equivalent locally uniformly convex norm, see \cite{DGZ}. Such a norm on the space $C(G)$, where $G$ is a compact group, has already been constructed by Aleksandrov in \cite{Alex}.


\section{Preliminaries}

We use standard notation concerning topology, set theory and Banach spaces. By a {\em map} we mean a continuous map. 
A {\em projection} in a Banach space $E$ is a bounded linear operator $\map PEE$ such that $PP=P$. One says that $F$ is {\em complemented} in $E$ if $F=\im P:=\setof{ Px}{x\in E}$ for some projection $\map PEE$. More precisely, $F$ is {\em $k$-complemented} if $\norm P\loe k$.
Let $T$ be a linear operator between subspaces of Banach spaces of the form $C(K)$. Then $T$ is called {\em regular} if $T$ is {\em positive}, i.e. $Tf\goe0$ whenever $f\goe0$, and $T1=1$, where $1$ denotes the constant function with value $1$ (so it is assumed that this function belongs to the domain of $T$).

Let $X,Y$ be two compact spaces and assume $\map fXY$ is a continuous surjection. We denote by $f^*$ the operator $\map S{C(Y)}{C(X)}$ defined by $S(\psi)= \psi f$, for $\psi\in C(Y)$. Clearly, $f^*$ is linear and provides an isometric embedding of $C(Y)$ into $C(X)$. One usually identifies $C(Y)$ with the subspace of $C(X)$, via the quotient map $f$.
A {\em regular averaging operator} associated with $f$, is a regular linear operator $\map T{C(X)}{C(Y)}$ satisfying $T(\psi f) = \psi$ for every $\psi\in C(Y)$.
Observe that given a regular averaging operator $\map T{C(X)}{C(Y)}$, the map $P=f^* T$ is a 
regular
(in particular: norm one) projection of $C(X)$ onto the subspace $\im f^* =\setof{ \psi f}{\psi\in C(Y)}$, isomorphic to $C(Y)$.
Regular averaging operators were introduced and studied by Pe\l czy\'nski \cite{Pelczynski}, motivated by Milyutin's Lemma \cite{Milyutin}, which says that there exists a continuous surjection of the Cantor set onto the unit interval admitting such an operator.

Given a compact space $K$, we denote by $P(K)$ the space of all regular probability measures on $K$. In other words, 
$$P(K)=\setof{\mu\in C(K)^*}{\norm\mu=1\text{ and }\mu(\phi)\goe0\text{ for every }\phi\goe0}.$$
We shall always consider $P(K)$ with the weak-star topology inherited from $C(K)^*$.
Every continuous map of compact spaces $\map fXY$ induces a map $\map{P(f)}{P(X)}{P(Y)}$ defined by $P(f)(\mu)(\phi)=\mu(\phi f)$, $\phi\in C(Y)$. By this way $P$ becomes a functor, usually called the {\em probability measures functor}. Note that $P(K)$ is a convex compact subset of $C(K)^*$ and it is second countable whenver $K$ is so. 
Given a second countable compact $K$, a special case of Michael's Selection Theorem says that every lower semi-continuous map $\Phi$, defined on a paracompact space $X$, whose values are closed convex subsets of $P(K)$, has a continuous selection, i.e. a map $\map hX{P(K)}$ such that $h(x)\in\Phi(x)$ for every $x\in X$.

A Banach space $E$ is {\em $k$-\pli} if there are a linearly dense set $X\subs E$ and a $k$-norming set $Y\subs E^*$ such that for every $y\in Y$ the set
$\setof{x\in X}{y(x)\ne0}$
is countable.
Recall that $Y$ is {\em $k$-norming} if
$\norm v \loe k \sup\setof{|y(v)|/\norm y}{y\in Y}$ for every $v\in E$.
We shall be interested in $1$-\pli~spaces.
In the case of density $\aleph_1$, $1$-\pli~spaces are characterized as Banach spaces with a {\em projectional resolution of the identity}, i.e. with a sequence $\sett{P_\al}{\alom}$ of norm one projections onto separable subspaces satisfying the following conditions:
\begin{enumerate}
	\item $\al<\beta\implies P_\al = P_\al P_\beta = P_\beta P_\al$;
	\item $\bigcup_{\alom}\im P_\al$ is the whole space and $\im P_\delta =\cl(\bigcup_{\al<\delta} \im P_{\al+1})$ for every limit ordinal $\delta$,
\end{enumerate}
where $\omega_1$ denotes the first uncountable ordinal.
For details we refer to Kalenda's survey \cite{Kalenda}. The well known notion of a projectional resolution of the identity is defined for an arbitrary Banach space $E$, where it is required that the density of $P_\al E$ does not exceed the cardinality of $\al+\omega$, see e.g. \cite{Fabian, DGZ}.

\pli~spaces are closely related to Valdivia compacta. Recall that a compact space $K$ is called {\em Valdivia compact} if $K\subs[0,1]^\kappa$ so that $K\cap \Sigma(\kappa)$ is dense in $K$, where $\Sigma(\kappa)$ is the {\em $\Sigma$-product} of $\kappa$ copies of $[0,1]$, i.e. $\Sigma(\kappa)=\setof{x\in[0,1]^\kappa}{|\setof{\al}{x(\al)\ne0}|\loe\aleph_0}$.
Let us recall that the Banach space $C(K)$ is $1$-\pli~whenever $K$ is Valdivia compact. On the other hand, straight from the definition it follows that for a $1$-\pli~Banach space $E$, the dual unit ball of $E$ endowed with the weak-star topology is Valdivia compact. For details we refer to \cite{Kalenda}.

We are going to use inverse sequences of compact spaces, so we briefly recall the necessary definitions.
Let $\delta$ be an infinite limit ordinal. An {\em inverse sequence} of length $\delta$ is a triple of the form $\S=\invsys Xp\delta\al\beta$, where for each $\al<\delta$, $X_\al$ is a topological (typically: compact) space and for each $\al<\beta<\delta$, $\map{p^\beta_\al}{X_\beta}{X_\al}$ is a continuous (typically: quotient) map, called a {\em bonding map}. Moreover, the following compatibility is required: $p^\gamma_\al = p^\beta_\al p^\gamma_\beta$ for every $\al<\beta<\gamma<\delta$.
The {\em limit} of $\S$ is a space $X=\liminv\S$ together with maps $\map{p_\al}X{X_\al}$ ($\al<\delta$) satisfying the following condition:
given a topological space $Y$ and a collection of maps $\sett{f_\al}{\al<\delta}$ such that $\map{f_\al}Y{X_\al}$ and $f_\al = p^\beta_\al f_\beta$ for every $\al<\beta<\delta$, there exists a unique map $\map fYX$ such that $p_\al f = f_\al$ holds for every $\al<\delta$.
The maps $p_\al$ are called {\em projections}.
Typically, $\liminv\S$ is represented as $X=\setof{x\in\prod_{\al<\delta}}{p^\beta_\al (x(\beta))=x(\al)\text{ for every }\al<\beta<\delta}$, where $p_\al$ is the projection onto $\al$-th coordinate. The inverse sequence $\S=\invsys Xp\delta\al\beta$ is {\em continuous} if for every limit ordinal $\gamma<\delta$ the space $X_\gamma$ together with the collection $\sett{p^\gamma_\al}{\al<\gamma}$ is the limit of the sequence $\invsys Xp\gamma \al\beta$.

We recall the definition of the class \R, introduced in \cite{BKT}. It is the smallest class of (compact) spaces that contains all metric compacta and which is stable under limits of continuous inverse sequences whose bonding maps are retractions. Every Valdivia compact has a decomposition into a continuous inverse sequence of retractions onto smaller Valdivia compacta (see e.g. \cite{Kalenda}), therefore it belongs to \R. For more information concerning class \R~and its properties, we refer to \cite{K_classR}. 

A {\em Dugundji compact} is a compact space $K$ which is an {\em absolute extensor} for the class of all $0$-dimensional compact spaces; that is: given a $0$-dimensional compact $X$ and a continuous map $\map fAK$ defined on a closed subset of $X$, there exists a continuous map $\map FXK$ such that $F\rest A=f$. We shall use the following useful characterization of Dugundji compacta, due to Haydon \cite{Haydon}: 

\begin{twH}
Let $K$ be a compact space. Then $K$ is Dugundji compact if and only if $K=\liminv\S$, where $\S=\invsys Kp\kappa\xi\eta$ is a continuous inverse sequence such that $K_0$ is metrizable and each $p^{\xi+1}_\xi$ is an open surjection with a metrizable kernel.
\end{twH}

Recall that a quotient map of compact spaces $\map fXY$ has a {\em metrizable kernel} if there is a map $\map hXZ$ such that $Z$ is second countable and the diagonal map $f\triangle h$ is one-to-one. Equivalently: there exists a second countable space $Z$ such that $X$ embeds into $Y\times Z$ so that $f$ is homeomorphic to the projection onto the first coordinate.

Interesting and important examples of Dugundji spaces are compact groups, see Uspenskij's article \cite{Uspenskij}. Let us note that $0$-dimensional Dugundji compacta are Valdivia \cite{KM}, although by \cite{KU} there exist compact Abelian groups which are not in the class \R.

The following lemma is well known. We give the proof for the sake of completeness.

\begin{lm}\label{w3jtrpqwjrp}
Assume $X,Y$ are compact spaces and $\map f{X}{Y}$ is an open surjection with a metrizable kernel. Then $f$ admits a regular averaging operator.
\end{lm}

\begin{pf}
Let $Q$ be a metric compact such that $X\subs Y\times Q$ and $f$ is the projection onto the first coordinate. Let $\map\pi{X}Q$ denote the projection onto the second coordinate. Define a multifunction $\Phi$, from $Y$ to the power set of $Q$, by setting
$$\Phi(y) = \img \pi{f^{-1}(y)}.$$
Then $\Phi$ has nonempty compact values. Since $f$ is open, $\Phi$ is lower semi-continuous.
We identify $Q$ with a suitable subset of the metrizable locally convex space $C(Q)^*$, endowed with the weak-star topology. By Michael's Selection Theorem, there exists a continuous map $\map {h_0}Y{C(Q)^*}$ such that $h_0(y)\in \cl_*(\conv \Phi(y))$ for every $y\in Y$, where $\cl_*$ denotes the weak-star closure. It follows that $h_0(y)$ is a probability measure whose support is contained in $\img\pi{f^{-1}(y)}$. Now define $h(y)\in P(X)$ by setting 
$$h(y)(\psi) = h_0(\psi_y),$$
where $\psi_y\in C(Q)$ is defined by $\psi_y(t)=\psi(y,t)$. The map $y\mapsto \psi_y$ is continuous with respect to the norm topology on $C(Q)$, therefore $\map h Y{P(X)}$ is continuous with respect to the weak-star topology on $P(X)$. 
Now define $\map T{C(X)}{C(Y)}$ by
$$(T\psi)(y) = \int_X \psi \, d h(y).$$
By the continuity of $h$, $T\psi$ is indeed a continuous function. Thus $T$ is a regular linear operator. Now assume $\psi = \phi f$. Then $\psi$ has constant value $\phi(y)$ on the set $f^{-1}(y)$. Recalling that the support of $h(y)$ is contained in $f^{-1}(y)$, we deduce that $(T\psi)(y) = \phi(y)$. Thus $T$ is a regular averaging operator.
\end{pf}

A collection of sets $\sett{S_\al}{\al<\lam}$ is a {\em chain} if $S_\al\subs S_\beta$ whenever $\al<\beta$. A chain $\sett{E_\al}{\al<\lam}$ of closed subspaces of a Banach space is {\em continuous} if $E_\delta=\cl(\bigcup_{\al<\delta}E_\al)$ for every limit ordinal $\delta<\lam$.

\begin{lm}[\cite{K2006}]\label{klocki}
Let $E$ be a Banach space and assume $\sett{E_\al}{\al<\lam}$ is a continuous increasing chain of closed subspaces of $E$ with $E=\cl(\bigcup_{\al<\lam}E_\al)$ ($\lam$ is a limit ordinal). Assume that for each $\al<\lam$, $\map{R_\al}{E_{\al+1}}{E_\al}$ is a norm one projection. Then there exists a sequence $\sett{P_\al}{\al<\lam}$ of projections of $E$ such that
\begin{enumerate}
	\item[(1)] $\norm{P_\al}=1$ and $P_\al E = E_\al$,
	\item[(2)] $\al\loe \beta<\lam\implies P_\al P_\beta = P_\beta P_\al = P_\al$.
	\item[(3)] $P_\al\rest E_{\al+1} = R_\al$.
\end{enumerate}
If, additionally, $E$ is of the form $C(K)$ and each $R_\al$ is regular, then we may assume that each $P_\al$ is regular.
\end{lm}

\begin{pf} We construct inductively norm one linear projections $\map{P^\beta_\al}{E_\beta}{E_\al}$, where $\al<\beta\loe\lam$ and $E_\lam=E$, satisfying the following condition
\begin{equation}
\al<\beta<\gamma\implies P^\beta_\al P^\gamma_\beta=P^\gamma_\al. 
\tag{*}\end{equation}
Suppose $\delta>0$ and $P^\beta_\al$ have been defined for all $\al\loe\beta<\delta$. If $\delta=\rho+1$ then define $P^{\rho+1}_\al = P^\rho_\al Q_\al$ for $\al<\delta$.
Assume now that $\delta$ is a limit ordinal. Let $D=\bigcup_{\xi<\delta}E_\xi$. Then $D$ is a dense linear subspace of $E_\delta$. Fix $\al<\delta$ and define
$P^\delta_\al(x)=P^\xi_\al(x)$, where $\xi$ is any ordinal satisfying $\al<\xi<\delta$ and $x\in E_\xi$.
Then ${P^\delta_\al}$ is a well defined norm one projection of $D$ onto $E_\al$, therefore it extends uniquely onto $E_\delta$. Given $\al<\beta<\delta$, the formula $P^\beta_\al P^{\delta}_\beta(x) = P^\delta_\al(x)$ is valid for every $x\in D$, therefore by continuity it holds for every $x\in E_\delta$. 

Now define $P_\al = P^\lam_\al$. Clearly (1) holds. Fix $\al<\beta$. Condition (*) for $\gamma=\lam$ says that $P_\al = P_\al^\beta P_\beta$, therefore $P_\al\rest E_\beta = P^\beta_\al$. Thus $P_\al P_\beta = P_\al$ and also (3) holds, because $P^{\al+1}_\al=R_\al$. On the other hand, $P_\beta P_\al = P_\al$, because $E_\al \subs E_\beta$. This shows (2).

Finally, in case where $E=C(K)$, it suffices to recall that regular operators are closed under compositions and pointwise limits.
\end{pf}

\section{Main result}


\begin{tw}\label{aewfacsfas}
Let $K$ be a Dugundji compact of weight $\aleph_1$. Then $C(K)$ has a projectional resolution of the identity $\sett{P_\al}{\alom}$ such that each $P_\al$ is regular. In particular, $C(K)$ is $1$-\pli~and the space of probability measures $P(K)$ is Valdivia compact.
\end{tw}

\begin{pf}
By Haydon's Theorem, $K=\liminv \invsys Kp{\omega_1}\al\beta$, where each $K_\al$ is a compact metric space, each $p^\beta_\al$ is an open surjection and the sequence is continuous. It suffices to show that for every $\al<\omega_1$, $C(K_\al)$ is complemented by a regular projection of $C(K_{\al+1})$, where we identify each $C(K_\al)$ with $\img{p_\al^*}{C(K_\al)}$ ($\map{p_\al}{K}{K_\al}$ is the projection). Indeed, Lemma \ref{klocki} will give us a projectional resolution of the identity $\sett{P_\al}{\alom}$ on $C(K)$ which consists of regular operators. Thus $C(K)$ is $1$-\pli~(see \cite{Kalenda}). Since $P_\al$'s are regular, the dual operators $P_\al^*$ provide a continuous inverse sequence of retractions of $P(K)$ onto metrizable subspaces which, by \cite[Corollary 4.3]{KM}, shows that $P(K)$ is Valdivia compact.

Fix $\alom$. Since $q:=q^{\al+1}_{\al}$ is open, by Lemma \ref{w3jtrpqwjrp}, it admits a regular averaging operator $\map T{C(K_{\al+1})}{C(K_{\al})}$. Then $R_\al := q^* T$ is a norm one regular projection of $C(K_{\al+1})$ onto $\img{q^*}{C(K_{\al})}$, identified with $C(K_{\al})$.
\end{pf}

The following statement gives a negative answer to a question of O. Kalenda \cite[Question 5.1.10(i)]{Kalenda}.

\begin{tw} There exists a compact space $K$ of weight $\aleph_1$, such that $K\notin \R$, while $P(K)$ is Valdivia compact and $C(K)$ is $1$-\pli.
\end{tw}

\begin{pf}
Let $K$ be the Abelian group described in \cite{KU}. Then $\w(K)=\aleph_1$, $K\notin\R$ and $K$ is a Dugundji space (being a compact group \cite{Uspenskij}). Thus, by Theorem \ref{aewfacsfas}, $P(K)$ is Valdivia and $C(K)$ is $1$-\pli.
\end{pf}

\begin{uwgi}
One should point out two things concerning Theorem \ref{aewfacsfas}. First, we did not have to show that $P(K)$ is Valdivia compact, because it follows from \cite[Theorem 5.1.2]{Kalenda}. Second, Theorem \ref{aewfacsfas} is valid (with the same proof) for every compact space $K$ which can be represented as $K=\liminv\S$, where $\S=\invsys Kp{\omega_1}\al\beta$ is a continuous inverse sequence of metric compacta such that for every $\alom$ the map $\map{p^{\al+1}_\al}{K_{\al+1}}{K_\al}$ admits a (not necessarily regular!) norm one averaging operator, i.e. such that $C(K_\al)$ is $1$-complemented in $C(K_{\al+1})$, under the suitable identification. This statement fails when $\omega_1$ is replaced by $\omega_2$ (the first ordinal of cardinality $\aleph_2$). For example, consider $K=\omega_2+1$ as the linearly ordered space. It has been proved by Kalenda \cite{Ondrej} that $C(K)$ is not a \pli~space. On the other hand, $K$ is the limit of a continuous inverse sequence of spaces of weight $\aleph_1$ in which all bonding maps are retractions (define $K_\al=\al+1$, $p^{\al+1}_\al\rest K_\al=\id_{K_\al}$ and $p^{\al+1}_\al(\al+1)=\al$).
\end{uwgi}

\begin{uwgi} It is well known that for every Dugundji compact $K$ there exists a continuous surjection $\map f{2^\kappa}K$ which admits a regular averaging operator. This gives a one-complemented isometric embedding of $C(K)$ into $C(2^\kappa)$ and a retraction of $P(2^\kappa)$ onto $P(K)$. The space $P(2^\kappa)$ is Valdivia compact and the space $C(2^\kappa)$ is $1$-\pli. 
\end{uwgi}


\end{document}